\documentclass[aos,MSNbibl,nameyear,dvips]{arximspdf}
\usepackage{mathrsfs}
\usepackage{graphicx}

\doi{10.1214/14-AOS1231} 
\volume{42}
\issue{5}
\pubyear{2014}
\firstpage{1725}
\lastpage{1750}
\docsubty{FLA}

\makeatletter
\def\mid{|}
\renewcommand{\citep}[1]{(\citeauthor{#1} \citeyear{#1})}
\newcommand{\eqref}[1]{(\ref{#1})}
\newtheorem{theorem}{Theorem}[section]
\newtheorem{lemma}[theorem]{Lemma}

\newtheorem{corollary}[theorem]{Corollary}
\makeatother

\begin{document}
\begin{frontmatter}

\title{A central limit theorem for general orthogonal array based space-filling designs}
\runtitle{A CLT for OA based space-filling designs}

\begin{aug}
\author[a]{\fnms{Xu}~\snm{He}\ead[label=e1]{hexu@amss.ac.cn}}
\and
\author[b]{\fnms{Peter Z.~G.}~\snm{Qian}\corref{}\thanksref{T1}\ead[label=e2]{peterq@stat.wisc.edu}}
\affiliation{Chinese Academy of Sciences and University of Wisconsin-Madison}

\thankstext{T1}{Supported by NSF Grant DMS-10-55214.}
\address[a]{Academy of Mathematics and Systems Science\\
Chinese Academy of Sciences\\
Zhongguancundonglu 55\\
Beijing 100190\\
China\\
\printead{e1}}
\address[b]{Department of Statistics\\
University of Wisconsin-Madison\\
1300 University Avenue\\
Madison, Wisconsin 53706\\
USA \\
\printead{e2}}

\runauthor{X. He and P. Z. G. Qian}
\end{aug}

\received{\smonth{12} \syear{2013}}
\revised{\smonth{4} \syear{2014}}

%
\begin{abstract}
Orthogonal array based space-filling designs
(Owen [\textit{Statist. Sinica} \textbf{2} (1992a) 439--452];
Tang [\textit{J. Amer. Statist. Assoc.} \textbf{88} (1993) 1392--1397])
have become popular in computer experiments, numerical integration,
stochastic optimization and uncertainty quantification.
As improvements of ordinary Latin hypercube designs,
these designs achieve stratification in multi-dimensions.
If the underlying orthogonal array has strength $t$, such designs achieve
uniformity up to $t$ dimensions.
Existing central limit theorems are limited to these designs
with only two-dimensional stratification based on strength two
orthogonal arrays.
We develop a new central limit theorem for these designs that possess
stratification in
arbitrary multi-dimensions associated with orthogonal arrays of general
strength.
This result is useful for building confidence statements for such
designs in various statistical
applications.
\end{abstract}

%

\end{frontmatter}

\section{Introduction}\label{sec1}

Latin hypercube designs achieve maximum uniformity in univariate
margins [\citet{McKay:Beckman:Conover:1979}].
Orthogonal arrays based Latin hypercube designs [\citet{Tang:1993}],
called \textit{U designs},
improve upon them by achieving uniformity in multivariate dimensions. %
Another type of orthogonal array based design is the randomized
orthogonal array [\citet{Patterson:1954,Owen:1992:ROA}].
The two classes of designs are widely used in computer experiments,
numerical integration, stochastic optimization and uncertainty quantification.

Consider a $K$-dimensional numerical integration problem
\[
\mu= E \bigl\{f(x)\bigr\} = \int_{[0,1)^K} f(x) \,dx.
\]
After evaluating $f$ at $N$ runs, $X_1,\ldots,X_N$, $\mu$ is estimated by
%
\begin{equation}
\hat\mu= N^{-1} \sum_{i=1}^N f
\bigl(X_i^1,\ldots,X_i^K\bigr),
\label{eqn:muhat}
\end{equation}
where $X_i^k$ is the $k$th dimension of $X_i$.
\citet{Tang:1993} gives a variance formula of $\hat\mu$ for a $U$ design,
and \citet{Owen:1994} derives variance formulas for a randomized
orthogonal array free of coincidence defect.
Methods to numerically estimate this variance are discussed in
\citeauthor{Owen:1992:ROA} (\citeyear{Owen:1992:ROA,Owen:1994}).

When an orthogonal array based space-filling design is used
in numerical integration, stochastic optimization [\citet
{Birge:2011,Shapiro:2009,Tang:Qian:2010}],
uncertainty quantification [\citet{Xiu:2010}] and other applications,
one is often interested in a central limit theorem for deriving a
confidence statement.
Derivation of a central limit theorem for such designs is a very
challenging problem
because of their complicated combinatorial structure and sophisticated
dependence
across the rows after randomization.
\citeauthor{Loh:1996} (\citeyear{Loh:1996,Loh:2008}) was first to address this problem
and derived central limit theorems for these designs associated with
orthogonal arrays of index one and strength two, which achieve
uniformity up to two-dimensional projections.
In \citet{Loh:2008}, the integrand is assumed to be Lipschitz
continuous mixed partial of order $K$.

Different from the work of \citeauthor{Loh:1996}
(\citeyear{Loh:1996,Loh:2008}), we
propose a new approach to construct a new central limit theorem for
orthogonal array based space-filling designs.
This approach works for these designs that
achieve uniformity in
arbitrary multi-dimensions associated with orthogonal arrays of general
strength.
As in \citet{Owen:1994}, we assume the underlying orthogonal array is
free of coincidence defect.
Let $\lambda$ and $n$ denote the index and the number of levels for the
orthogonal array, respectively. As $N$ tends to infinity, we assume
$\lambda$ is fixed or $\lambda/n$ tends to zero.
Our method is inspired by the method of moments used in \citet
{Owen:1992:CLT} for ordinary Latin hypercube designs but with new
combinational techniques to deal with the complexity of orthogonal arrays.
Section~\ref{sec:def} presents useful definitions and notation.
Sections~\ref{sec:ROA} and \ref{sec:UD} provide central limit theorems
for orthogonal array based space-filling designs. Section~\ref{sec:simu} gives numerical illustration of the derived theoretical results.
Section~\ref{sec:con} concludes with some brief discussion.

\section{Definitions and notation}\label{sec:def}

An $N$ by $K$ matrix is said to be a Latin hypercube if each of its
columns consists of $\{0,1,\ldots,N-1\}$.
A uniform permutation on a set of $a$ numbers is randomly generated
with all $a!$ permutations equally probable.
An ordinary Latin hypercube design~[\citet{McKay:Beckman:Conover:1979}]
is constructed by
\[
X_i^k = \pi_k(i) / N + \eta_i^k
/ N, \label{eqn:const:LH}
\]
where the $\pi_k$ are uniform permutations on $\{0,1,\ldots,N-1\}$,
the $\eta_i^k$ are generated from uniform distributions on $[0,1)$
and the $\pi_k$ and the $\eta_i^k$ are generated independently.

An $N$ by $K$ matrix is said to be an orthogonal array $OA(N,K,n,h)$ if
its entries are from $0,1,\ldots,n-1$
and for any $p \leq h$ columns of the matrix, the $n^p$ combinations of
values appear exactly the same number of times in rows~[\citet
{Hedayat:Sloane:Stufken:1999}].
For an $OA(N,K,n,h)$, if additionally no two rows from any $N \times
(h+1)$ submatrices are the same, the orthogonal array is said to be
free of coincidence defect [\citet{Owen:1994}].
For illustration, Table~\ref{tab:OA18} gives an $OA(18,6,3,2)$ of index
two and free of coincidence defect.

\begin{table}
\caption{An orthogonal array with 18 runs}\label{tab:OA18}
\begin{tabular*}{150pt}{@{\extracolsep{\fill}}lccccc@{}}
\hline
0&0&0&0&0&0 \\
1&1&1&1&1&1 \\
2&2&2&2&2&2 \\
0&0&1&2&1&2 \\
1&1&2&0&2&0 \\
2&2&0&1&0&1 \\
0&1&0&2&2&1 \\
1&2&1&0&0&2 \\
2&0&2&1&1&0 \\
0&2&2&0&1&1 \\
1&0&0&1&2&2 \\
2&1&1&2&0&0 \\
0&1&2&1&0&2 \\
1&2&0&2&1&0 \\
2&0&1&0&2&1 \\
0&2&1&1&2&0 \\
1&0&2&2&0&1 \\
2&1&0&0&1&2\\
\hline
\end{tabular*}
\end{table}

Let $H$ denote an $OA(N,K,n,h)$ with the $(i,k)$th element $H_i^k$.
A randomized orthogonal array~[\citet{Owen:1992:ROA}] based on $H$ is
constructed by
%
\begin{equation}
X_i^k = \pi_k\bigl(H_{\gamma^{-1}(i)}^k
\bigr) / n + \eta_i^k / n, \label{eqn:const:ROA}
\end{equation}
where
the $\gamma$ is a uniform permutation on $\{1,\ldots,N\}$,
the $\pi_k$ are uniform permutations on $\{0,1,\ldots,n-1\}$,
the $\eta_i^k$ are generated from the uniform distribution on $[0,1)$
and the $\gamma$, the $\pi_k$ and the $\eta_i^k$ are generated independently.

Compared with (\ref{eqn:const:ROA}), a $U$ design~[\citet{Tang:1993}] based
on $H$ is constructed with one additional step,
%
\begin{equation}
X_i^k = \pi_k\bigl(H_{\gamma^{-1}(i)}^k
\bigr) / n + \alpha_{\gamma^{-1}(i)}^k / N + \eta_i^k
/ N, \label{eqn:const:UD}
\end{equation}
where
the $\gamma$ is a uniform permutation on $\{1,\ldots,N\}$,
the $\pi_k$ are uniform permutations on $\{0,1,\ldots,n-1\}$,
all the $\alpha_{i}^k$'s related to entries in the $k$th column with
level $x$ in $H$ consist of a permutation of $\{0,1,\ldots,N/n-1\}$,
the $\eta_i^k$ are generated from uniform distributions on $[0,1)$
and the $\gamma$, the $\pi_k$, the $\alpha_{k,x}^i$ and the $\eta_i^k$
are generated independently.

\begin{figure}

\includegraphics{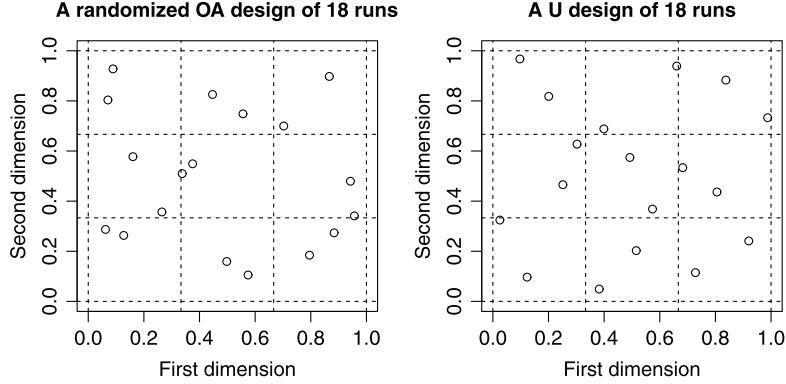}

\caption{Bivariate projections to the first two dimensions of a
randomized orthogonal array design and a $U$ design generated from the
orthogonal array in Table~\protect\ref{tab:OA18}.
For both designs, each of the nine squares by dashed lines contains
exactly two points.
Furthermore, for the $U$ design, each of the 18 equally spaced intervals
of $[0,1)$ contains exactly one point.}
\label{fig:exp18}
\end{figure}

For illustration, let $H$ be the orthogonal array in Table~\ref{tab:OA18}.
We generate a randomized orthogonal array and a $U$ design based on $H$.
The bivariate projections to the first two dimensions of the two
designs are depicted in Figure~\ref{fig:exp18}.
For both designs, each of the nine squares by dashed lines contains
exactly two points.
Furthermore, for the $U$ design, each of the 18 equally spaced intervals
of $[0,1)$ contains exactly one point.

Next, we introduce the functional analysis of variance decomposition
[\citet{Owen:1994}].
Let $F$ be the uniform measure on $[0,1)^K$ with $dF=\break \prod_{k=1}^K
\,dF_{\{k\}}$,
where $F_{\{k\}}$ is the uniform measure on $[0,1)$.
Under the assumption $f$ is a continuous function in $[0,1]^K$,
$f$ is bounded and has finite variance $\int f(x)^2 \,dF$.
Express $f$ as
\[
f(x)= \mu+ \sum_{\phi\subset u \subseteq\{1,\ldots,K\}} f_u(x),
\]
where
$ \mu= \int f(x) \,dF $
and $f_u$ is defined recursively via
\[
f_u(x) = \int\biggl\{ f(x) - \sum_{v \subset u}
f_v(x) \biggr\} \,dF_{\{1,\ldots,K\}
\setminus u}.
\]
If $u \cap v \neq\phi$,
%
\begin{equation}
\int_v f_u \,dx = 0. \label{eqn:int:0}
\end{equation}

Following \citet{Owen:1994},
for the two classes of designs of strength $h$ without coincidence defect,
the $f_u$ part with $|u|\leq h$ is balanced out from the design.
The remaining part $r$ of $f$ is defined via
%
\begin{equation}
f(x)= \mu+ \sum_{0<|u|\leq h} f_u(x) +r.
\label{eqn:int:r}
\end{equation}
The variance of $\hat\mu$ from (\ref{eqn:muhat}) is
\[
\mbox{var}(\hat\mu)= N^{-1} \int r(X)^2 \,dF(X) + o
\bigl(N^{-1}\bigr).
\]

Let $I(\cdot)$ be the indicator function.
For a real number $x$, let $\lfloor x \rfloor$ be the largest integer
no greater than $x$,
and the subdivision of $x$ with length $1/z$ is
%
\begin{equation}
\delta_z(x) = \bigl[ \lfloor z x \rfloor/z, \bigl(\lfloor z x \rfloor+1
\bigr) /z \bigr). \label{eqn:delta}
\end{equation}
Let $|D|$ be the volume of region $D$.
Let $E_\mathrm{ IID}$, $E_\mathrm{ ROA}$ and $E_\mathrm{ UD}$ be the expectation of
a function from samples generated identically and independently, from a
randomized orthogonal array and from a $U$ design, respectively.

\section{A central limit theorem for randomized orthogonal arrays}\label
{sec:ROA}

We now derive a central limit theorem for randomized orthogonal arrays.
Assume $f$ is a continuous function from $[0,1]^K$ to $\mathscr{R}$.
Let $H$ be an $OA(N,K,n,h)$ free of coincidence defect and $\lambda= N/n^h$.
Take $X_1,\ldots,X_N$ in (\ref{eqn:muhat}) to be the design points from
a randomized orthogonal array constructed in (\ref{eqn:const:ROA}).
For fixed $K$ and $h$, we suppose there is a sequence of $H$ such that
$N$ and $n$ tend to infinity with $\lambda/n$ tending to zero.
Lemma~\ref{lem:mom} on the method of moments~[\citet{Durrett:2010}] is
used throughout.

\begin{lemma}\label{lem:mom}
Suppose that $A_1,A_2,\ldots$ are random variables, and their
distribution functions $F_1,F_2,\ldots$ have finite moments.
Namely, for any $p=1,2,\ldots$ and $n=1,2,\ldots,$
\[
m_n^{(p)} = \int_{-\infty}^{+\infty}
x^p \,dF_n
\]
is finite. Suppose that $F$ is a distribution function with finite
moments. Namely,
\[
m^{(p)} = \int_{-\infty}^{+\infty} x^p
\,dF
\]
is finite.
Also assume
\[
\limsup_{p\to\infty} \bigl\{ \bigl(m^{(2p)}
\bigr)^{1/2p} / (2p) \bigr\} < \infty.
\]
Finally, suppose for any $p=1,2,\ldots,$
\[
\lim_{n \to\infty} m_n^{(p)} = m^{(p)}.
\]
Then $A_n$ converges in distribution to $F$.
\end{lemma}

Because the density function of multiple points among $X_1,\ldots,X_N$
is complicated,
we consider the conditional density $g=g(d_1,\ldots,d_K)$ of $X_s$
given other points $X_1,\ldots,X_{s-1}$, $s=1,\ldots,N$.
Unfortunately, the conditional density is not uniquely determined by
the definition of orthogonal arrays and $N,K,n,h$
and depends on the specific construction algorithm of $H$.
A key to overcome this difficulty is to express $g$ in big O terms.
Let $M_{s-1}$ denote an $(s-1)\times K$ matrix with the $(i,k)$th
element being $z$ if $z<i$, and $z$ is the smallest number such that
$X_i^k$ matches $X_z^k$, that is, $\lfloor n X_{i,k}\rfloor= \lfloor n
X_{z,k}\rfloor$.
If $X_i^k$ does not match to any other point $X_z^k$ with $z<i$, the
$(i,k)$th element of $M_{s-1}$ is defined to be zero and the first row
of $M_{s-1}$ is zero.
According to this definition, $M_{s-1}$ contains full information on
pairwise coincidence among $X_1,\ldots,X_{s-1}$.

\begin{lemma}\label{lem:ROA:g}
For a randomized orthogonal array in (\ref{eqn:const:ROA}), the
conditional density of $X_s$ given $X_1,\ldots,X_{s-1}$ is
%
\begin{eqnarray}\label{eqn:g}
&&g_s(d_1,\ldots,d_K)
\nonumber
\\[-8pt]
\\[-8pt]
\nonumber
&&\qquad= \sum
_{i_1,\ldots,i_K=0}^{s-1} b_s(i_1,
\ldots,i_K,M_{s-1}) I\bigl(d_1\in
D_{i_1}^1,\ldots,d_K\in D_{i_K}^K
\bigr),
\end{eqnarray}
where
$D_i^k = \delta_n(X_i^k)$ for $i=1,\ldots,s-1$ and $k=1,\ldots,K$,
$D_0^k=[0,1) \setminus\break \{ \bigcup_{i=1}^{s-1} \delta_n(X_i^k) \} $ for
$k=1,\ldots,K$ and
$b_s(\cdot)$ is a deterministic function on $d_1,\ldots, d_K, M_{s-1}$ with
\[
b_s(i_1,\ldots,i_k, M_{s-1}) = \cases{
1+O\bigl(n^{-1}\bigr), &\quad $|w|< h,$
\vspace*{2pt}\cr
O(1), &\quad $|w|= h,$
\vspace*{2pt}\cr
0, & \quad$|w|>h, \max\bigl(|w_1|,\ldots,|w_{s-1}|\bigr) >h,$
\vspace*{2pt}\cr
O\bigl(n^{|w|}/N\bigr), &\quad $\mbox{otherwise,}$}
\]
where $w(i_1,\ldots,i_K)$ is the dimensions of nonzero elements in
$(i_1,\ldots,i_K)$, $w=\{k\dvtx i_k\neq0\}$,
$w_z(i_1,\ldots,i_K) = \{k\dvtx i_k=z\}$ and $w=\bigcup w_z$.
\end{lemma}

\begin{figure}

\includegraphics{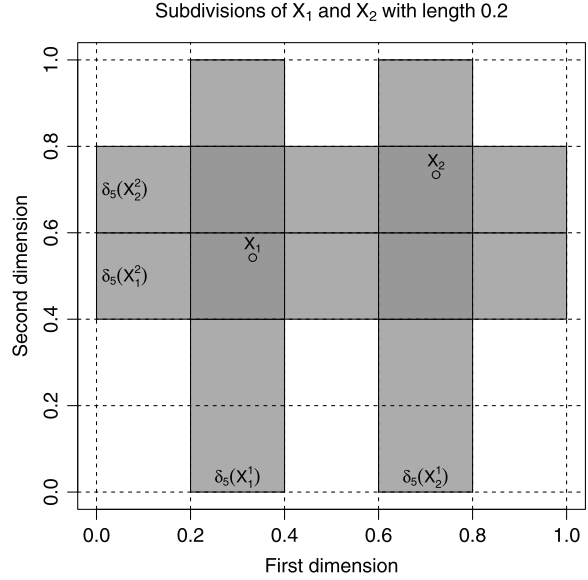}

\caption{The subdivisions of $X_1$ and $X_2$ with length $1/5$ for an
example with $n=5$, $h=2$, $K=2$ and $s=3$.
The white, light gray and gray regions represent the cases with
$|w|=0$, $|w|=1$ and $|w|=2$, respectively.}
\label{fig:D}
\end{figure}

Lemma~\ref{lem:ROA:g} shows that the conditional density is a constant
except in the subdivisions of $X_1,\ldots,X_{s-1}$
and $|w|$ indicates the number of dimensions that $X_s$ is inside the
subdivisions of any length.
For illustration, Figure~\ref{fig:D} displays subdivisions of $\delta
_n(X_i^k)$ for $n=5$, $h=2$, $K=2$ and $s=3$.
In this example, $X_1=(0\cdot332,0\cdot542)$ and $X_2=(0\cdot
722,0\cdot734)$.
The subdivisions of $X_1$ and $X_2$ are $\delta_5(X_1^1)=[0\cdot
2,0\cdot4)$, $\delta_5(X_1^2)=[0\cdot4,0\cdot6)$, $\delta
_5(X_2^1)=[0\cdot6,0\cdot8)$ and $\delta_5(X_2^2)=[0\cdot6,0\cdot8)$.
The regions with $|w|=0$, $|w|=1$ and $|w|=2$ are in white, light gray
and gray colors, respectively.
The proof of Lemma~\ref{lem:ROA:g} is given in the \hyperref[app]{Appendix}.

Next, we state two lemmas for the conditional expectation of $f(X_s)$
given points $X_1,\ldots,X_{s-1}$ from a randomized orthogonal array.
These lemmas parallel the results for ordinary Latin hypercube designs
in~\citet{Owen:1992:CLT} but use more complicated arguments.

\begin{lemma}\label{lem:ROA:f}
For any bounded function $f$ and $s >1$, as $N \to\infty$,
\[
E_\mathrm{ ROA}\bigl\{ f(X_s) \mid X_1,
\ldots,X_{s-1} \bigr\} = E_\mathrm{ IID} \bigl\{ f(X_s)
\bigr\} + O\bigl(n^{-1}\bigr).
\]
\end{lemma}

\begin{lemma}\label{lem:ROA:r}
Let
\[
\bar R= N^{-1} \sum_{i=1}^N
r(X_i).
\]
Then for any positive integer $p$,
\[
E_\mathrm{ ROA} \bigl\{\bigl(N^{1/2} \bar R\bigr)^p\bigr\}
= E_\mathrm{ IID} \bigl\{\bigl(N^{1/2} \bar R\bigr)^p\bigr
\} + o(1).
\]
\end{lemma}

Lemma~\ref{lem:ROA:f} is a direct consequence of Lemma~\ref{lem:ROA:g}.
The proof of Lemma~\ref{lem:ROA:r} is given in the \hyperref[app]{Appendix}.

We are now ready for our main theorem for randomized orthogonal arrays.

\begin{theorem}\label{thm:ROA}
Suppose that $f$ is a continuous function from $[0,1]^K$ to $\mathscr{R}$,
$\hat\mu$ in (\ref{eqn:muhat}) is based on a randomized orthogonal
array in (\ref{eqn:const:ROA}) without coincidence defect,
$\lambda$ is fixed or $\lambda=o(n)$.
Then, as $N \to\infty$,
\[
N^{1/2} ( \hat\mu- \mu) \to N \biggl(0, \int r(x)^2 \,dx
\biggr).
\]
\end{theorem}

\begin{pf}
The mean of $N^{1/2} ( \hat\mu- \mu)$ is 0 and the variance of
$N^{1/2} ( \hat\mu- \mu)$ tends to $\int r(x)^2 \,dx$.
From Lemma~\ref{lem:ROA:r}, for $p=1,2,\ldots,$
\[
E_\mathrm{ ROA} \bigl\{\bigl(N^{1/2} \bar R\bigr)^p\bigr\}
= E_\mathrm{ IID} \bigl\{\bigl(N^{1/2} \bar R\bigr)^p\bigr
\} + o(1).
\]
When the points are generated identically and independently, $N^{1/2}
\bar R$ follows a normal distribution with mean zero and variance
$\sigma^2=\int r(x)^2 \,dx$.
From~\citet{Owen:1980},
\[
E_\mathrm{ IID} \bigl\{\bigl(N^{1/2} \bar R\bigr)^p\bigr\}
= \cases{ %
0, &\quad $p=1,3,5,\ldots,$
\cr
\sigma^{p}(p-1)!!,
& \quad $p=2,4,6,\ldots.$ }
\]
Note that
\[
\limsup_{p\to\infty} \bigl(\sigma^p (p-1)!!
\bigr)^{1/p} / p =0.
\]
From Lemma~\ref{lem:mom}, $N^{1/2} \bar R$ from randomized orthogonal
array has the same limiting distribution as $N^{1/2} \bar R$ where the
points are generated identically and independently, which yields a
normal distribution.
\end{pf}

We can easily extend Theorem~\ref{thm:ROA} to a multivariate function
$f=(f_1,\ldots,\break f_P)$.
Parallel to (\ref{eqn:int:r}), define $r_i$ via
\[
f_i(x)= \mu_i + \sum_{0<|u|\leq h}
f_{i,u}(x) +r_i.
\]
The following theorem gives a central limit theorem for a multivariate $f$.

\begin{corollary}\label{thm:multi:ROA}
Suppose that $f$ is a continuous function from $[0,1]^K$ to $\mathscr{R}^P$,
$\hat\mu$ in (\ref{eqn:muhat}) is based on a randomized orthogonal
array in (\ref{eqn:const:ROA}) without coincidence defect,
$\lambda$ is fixed or $\lambda=o(n)$.
Then, as $N \to\infty$,
\[
N^{1/2} ( \hat\mu- \mu) \to N (0, \Sigma),
\]
where $\Sigma$ is a $P\times P$ matrix with the $(i,j)$th element
$\Sigma_{i,j} = \int r_i(x) r_j(x) \,dx $.
\end{corollary}

The normality of multivariate $f$ follows from the fact that any linear
combinations of $(f_1,\ldots,f_P)$ has a limiting normal distribution.

\section{A central limit theorem for $U$ designs}\label{sec:UD}

Next, we derive a central limit theorem for $U$ designs.
As before, we assume $f$ is a continuous function from $[0,1]^K$ to
$\mathscr{R}$.
Let $H$ be an $OA(N,K,n,h)$ free of coincidence defect and $\lambda= N/n^h$.
Take $X_1,\ldots,X_N$ in (\ref{eqn:muhat}) to be the design points from
a $U$ design constructed in (\ref{eqn:const:UD}).
For fixed $K$ and $h$, we suppose there is a sequence of $H$ such that
$N$ and $n$ tend to infinity with $\lambda/n$ tending to zero.
Analogous to Lemma~\ref{lem:ROA:g}, we first derive the conditional
density function of $X_s$ given $X_1,\ldots,X_{s-1}$.

\begin{lemma}\label{lem:UD:g}
For a $U$ design in (\ref{eqn:const:UD}) from $H$, the conditional
density of $X_s$ given $X_1,\ldots,X_{s-1}$ is
\[
g_s(d_1,\ldots,d_K) = \sum
_{i_1,\ldots,i_K=0}^{s-1} b_s(i_1,
\ldots,i_K,M_{s-1}) I\bigl(d_1\in
D_{i_1}^1,\ldots,d_K\in D_{i_K}^K
\bigr),
\]
where
$D_i^k = \delta_n(X_i^k) \setminus\{ \bigcup_{j=1}^{s-1} \delta_N(X_j^k)
\} $ for $i=1,\ldots,s-1$ and $k=1,\ldots,K$,
$D_i^k = \delta_N(X_{i-(s-1)}^k)$ for $i=s,\ldots,2s-2$ and $k=1,\ldots,K$,
$D_0^k=[0,1) \setminus\break \{ \bigcup_{j=1}^{s-1} \delta_n(X_j^k) \} $ for
$k=1,\ldots,K$ and
$b_s(\cdot)$ is a deterministic function on $d_1,\ldots,d_K, M_{s-1}$ with
\[
b_s(i_1,\ldots,i_k, M_{s-1}) = \cases{
0, & \quad$\mbox{there is a $k$ such that }
i_k> s-1,$
\vspace*{2pt}\cr
1+O\bigl(n^{-1}\bigr), &\quad $i_1,\ldots,i_K\leq
s-1, |w|< h,$
\vspace*{2pt}\cr
O(1), &\quad $i_1,\ldots,i_K\leq s-1, |w|= h,$
\vspace*{2pt}\cr
0, & \quad$|w|>h, \max\bigl(|w_1|,\ldots,|w_{s-1}|\bigr) >h,$
\vspace*{2pt}\cr
O\bigl(n^{|w|}/N\bigr), &\quad $\mbox{otherwise,}$}
\]
where $w(i_1,\ldots,i_K)$ is the dimensions of nonzero elements in
$(i_1,\ldots,i_K)$, $w=\{k\dvtx i_k\neq0\}$,
$w_z(i_1,\ldots,i_K) = \{k\dvtx i_k=z\}$ and $w=\bigcup w_z$.
\end{lemma}

The proof of Lemma~\ref{lem:UD:g} is given in the \hyperref[app]{Appendix}.
Analogous to Lemmas~\ref{lem:ROA:f} and~\ref{lem:ROA:r}, we state two
lemmas for the conditional expectation of $f(X_s)$ given points
$X_1,\ldots,X_{s-1}$ from a $U$ design.

\begin{lemma}\label{lem:UD:f}
For any bounded function $f$ and $s >1$, as $N \to\infty$,
\[
E_\mathrm{ UD}\bigl\{ f(X_s) \mid X_1,
\ldots,X_{s-1} \bigr\} = E_\mathrm{ IID} \bigl\{ f(X_s)
\bigr\} + O\bigl(n^{-1}\bigr).
\]
\end{lemma}

\begin{lemma}\label{lem:UD:r}
Let
\[
\bar R= N^{-1} \sum_{i=1}^N
r(X_i).
\]
Then for any positive integer $p$,
\[
E_\mathrm{ UD} \bigl\{\bigl(N^{1/2} \bar R\bigr)^p\bigr\}
= E_\mathrm{ IID} \bigl\{\bigl(N^{1/2} \bar R\bigr)^p\bigr
\} + o(1).
\]
\end{lemma}

Lemma~\ref{lem:UD:f} is a direct consequence of Lemma~\ref{lem:UD:g}.
A sketch to prove Lemma~\ref{lem:UD:r} is given in the \hyperref[app]{Appendix}.
A central limit theorem for $U$ designs is given below.

\begin{theorem}\label{thm:UD}
Suppose that $f$ is a continuous function from $[0,1]^K$ to $\mathscr{R}$,
$\hat\mu$ in (\ref{eqn:muhat}) is based on a $U$ design in (\ref
{eqn:const:UD}) without coincidence defect,
$\lambda$ is fixed or $\lambda=o(n)$.
Then, as $N \to\infty$,
\[
N^{1/2} ( \hat\mu- \mu) \to N \biggl(0, \int r(x)^2 \,dx
\biggr).
\]
\end{theorem}

\begin{pf}
$E\{ N^{1/2} ( \hat\mu- \mu) \} = 0$ and $\mbox{var}\{N^{1/2} ( \hat
\mu- \mu)\}$ tends to $\int r(x)^2 \,dx$.
From Lemma~\ref{lem:UD:r} and~\citet{Owen:1980}, for $p=1,2,\ldots,$
\begin{eqnarray*}
E_\mathrm{ UD} \bigl\{\bigl(N^{1/2} \bar R\bigr)^p\bigr\}
&=& E_\mathrm{ IID} \bigl\{\bigl(N^{1/2} \bar R\bigr)^p\bigr
\} + o(1) \\
&=& \cases{ %
 0+o(1), & \quad $p=1,3,5,\ldots,$
\vspace*{2pt}\cr
\sigma^{p}(p-1)!! +o(1), & \quad $p=2,4,6,\ldots,$ }
\end{eqnarray*}
where $\sigma^2 = \int r(x)^2 \,dx$ with
\[
\limsup_{p\to\infty} \bigl(\sigma^p (p-1)!!
\bigr)^{1/p} / p =0.
\]
From Lemma~\ref{lem:mom}, $N^{1/2} \bar R$ from $U$ design has the same
limiting distribution as $N^{1/2} \bar R$ where the points are
generated identically and independently, which yields a normal distribution.
\end{pf}

Similarly, the result can be extended to a multivariate $f$.

\begin{corollary}\label{thm:multi:UD}
Suppose that $f$ is a continuous function from $[0,1]^K$ to $\mathscr{R}^P$,
$\hat\mu$ in (\ref{eqn:muhat}) is based on a $U$ design in (\ref
{eqn:const:UD}) without coincidence defect,
$\lambda$ is fixed or $\lambda=o(n)$.
Then, as $N \to\infty$,
\[
N^{1/2} ( \hat\mu- \mu) \to N (0, \Sigma),
\]
where $\Sigma$ is a $P\times P$ matrix with the $(i,j)$th element
$\Sigma_{i,j} = \int r_i(x) r_j(x) \,dx $.
\end{corollary}

\section{Numerical illustration}\label{sec:simu}

\begin{figure}

\includegraphics{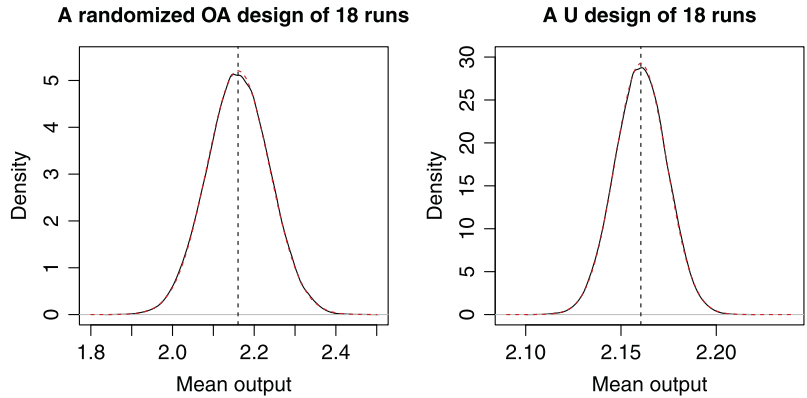}

\caption{Density plots of $\hat\mu$ based on a randomized orthogonal array
(left) and a $U$ design (right) from the orthogonal array given in
Table~\protect\ref{tab:OA18}, both of which are close to a normal distribution.}
\label{fig:CLT:18}
\end{figure}

We provide two numerical examples to validate the central limit
theorems in Sections~\ref{sec:ROA} and~\ref{sec:UD}.
In the first experiment, the orthogonal array with 18 runs, three
levels and strength two in Table~\ref{tab:OA18} of Section~\ref{sec:def} is used to generate a randomized orthogonal array and a $U$ design.
Consider estimating the mean output of a function [\citet{Cox:2001}]
\[
f = x_1/ \Bigl[ 2 \Bigl\{ \sqrt{1+\bigl(x_2+x_3^2
\bigr)x_4/x_1^2} -1 \Bigr\} \Bigr]
+x_1 +3x_4,
\]
where $x_1,\ldots,x_4$ follow the uniform distribution on $[0,1)$.
The true value of $\mu$ is approximately $2\cdot 160$, computed from a
large ordinary Latin hypercube design.
We compute $ \hat\mu= \sum_{i=1}^{18} f(X_i) /18 $ as in (\ref
{eqn:muhat}) for the two designs.
This procedure is repeated for 100,000 times.
The density plots of $\hat\mu$ for the two designs are shown in
Figure~\ref{fig:CLT:18},
where both distributions are close to a normal distribution.

\begin{figure}

\includegraphics{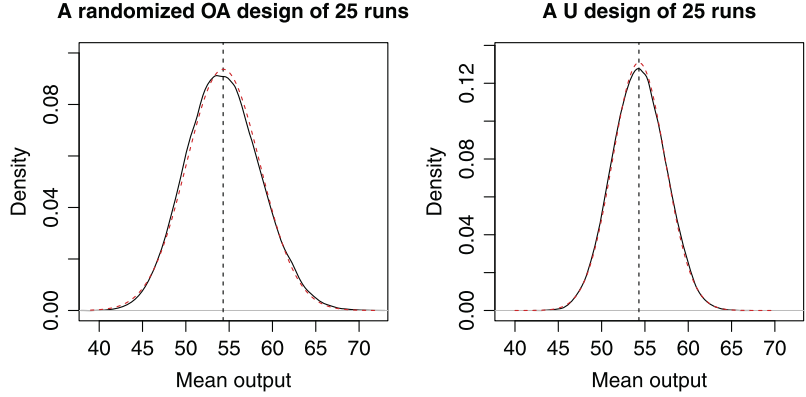}

\caption{Density plots of $\hat\mu$ based on a randomized orthogonal array
(left) and a $U$ design (right) from an orthogonal array with 25 runs,
both of which are close to a normal distribution.}
\label{fig:CLT:25}
\end{figure}

In the second experiment, an orthogonal array with 25 runs, five levels
and strength two is used for generating
a randomized orthogonal array and a $U$ design. We estimate the mean
output $\mu$ of the Branin function~[\citet{Branin:1972}]
\[
f = \biggl( x_2 - \frac{5.1}{4\pi^2}x_1^2 +
\frac{5}{\pi}x_1 -6 \biggr) + 10 \biggl(1-\frac{1}{8\pi}
\biggr) \cos(x_1) + 10
\]
on the domain $[-5,10] \times[0,15]$. The true value of $\mu$ is
approximately $54\cdot 31$, computed from a large grid design.
We compute $\hat\mu= \sum_{i=1}^{25} f(X_i) /25$
for the two designs. This procedure is repeated for 100,000 times.\vadjust{\goodbreak}
The density plots of $\hat\mu$ from the two designs are shown in
Figure~\ref{fig:CLT:25}, both of which are close to a normal distribution.

\section{Conclusions}\label{sec:con}
A new central limit theorem has been derived for orthogonal array based
space-filling designs. One might be interested in extending our
technique to derive a central limit theorem for scrambled nets [\citet
{Owen:1997}]. Another possible direction for future research is to use
this new result to study validation of sample average approximation
solutions for a stochastic program [\citet{Shapiro:2009}].
Finally, it is an important problem to estimate the variance $\int
r(x)^2\,dx$ from a $U$ design.

\begin{appendix}
\section*{Appendix}\label{app}

\subsection{Proof of Lemma~\texorpdfstring{\protect\ref{lem:ROA:g}}{3.2}} \label{sec:proof:ROA:g}
We first work on the $g_s(d_1,\ldots,d_K)$ on the cells $D=(D^1,\ldots,D^K)$,
where $D^k \in\{ [0,1/n), [1/n,2/n), \ldots, [(n-1)/n,1) \}$ for
$k=1,\ldots,K$.
Consider the matrix $\tilde H$ obtained by dropping rows $\gamma
^{-1}(1),\ldots,\break \gamma^{-1}(s-1)$ of $H$.
$g_s(d_1,\ldots,d_K)$ is nonzero if $(\lfloor n d_1 \rfloor,\ldots,\lfloor n d_K \rfloor)$ can be obtained from a row of $\tilde H$ by
some operators $\pi_k$,
which means $ \lfloor n d_1 \rfloor= \pi_1 (H_r^1),\ldots, \lfloor n
d_K \rfloor= \pi_K (H_r^K)$ for a row $H_r=(H_r^1,\ldots,H_r^K)$ in
$\tilde H$.
Let $x$ be the number of rows in $\tilde H$ from which $(\lfloor n d_1
\rfloor,\ldots,\lfloor n d_K \rfloor)$ can be obtained.
The value of $g_s(d_1,\ldots,d_K)$ is closely related to $x$ because
$X_s$ has the same probability $1/(N-(s-1))$ being permuted from each
row of $\tilde H$.

Because level permutations do not affect the result on whether two rows
of $H$ take same value in a particular column, $x$ is closely related
to $M_s$ and $w$.
Below we compute $x$ by types of $w$.

For the type of $|w|=0$, since there are at most $(s-1)(N/n-1)$ rows
taking value in $\bigcup_{i=1}^{s-1} \{ H_{\gamma^{-1}(i)}^k \}$ in the
$k$th column for $k=1,\ldots,K$,
$N - (s-1) - K(s-1)(N/n-1) \leq x \leq N-(s-1)$ and $x=N(1-O(n^{-1}))$.
Since the volume of cells for $w=\phi$ is $1-O(n^{-1})$ and $g_s(d)$ is
the same in such cells, $g_s(d_1,\ldots,d_K)=1+O(n^{-1})$.

For the type of $|w|=1$, without loss of generality, assume $w=\{1\}$
and $w_1=\{1\}$.
There are at least $N/n - (s-1)$ rows and at most $N/n -1$ rows taking
value $\{ H_{\gamma^{-1}(1)}^1 \}$ in the first column.
Out of those rows, there are at most $(s-1)(N/n^2-1)$ rows taking value
in $\bigcup_{i=1}^{s-1} \{ H_{\gamma^{-1}(i)}^k \}$ in the $k$th column
for $k=2,\ldots,K$.
Therefore, $x=N/n(1-O(n^{-1}))$.
Since the volume of cells for $w=\{1\}$ is $n^{-1}(1-O(n^{-1}))$ and
$g_s(d)$ is the same in such cells, $g_s(d_1,\ldots,d_K)=1+O(n^{-1})$.
Similarly, we obtain $g_s(d_1,\ldots,d_K)=1+O(n^{-1})$ for any $w$ with $|w|<h$.

For the type of $|w|=h$, there are at most $N/n^h$ rows in $\tilde H$
that match $X_s$ in $w$.
Since the volume of cells is $n^{-h}(1+O(n^{-1}))$, $g_s(d_1,\ldots,d_K)=O(1)$.

For the type of $|w|>h$, because $H$ is free of coincidence defect,
there is zero or one row in $\tilde H$ that matches $X_s$ in $w$.
Since the volume of cells is $n^{-|w|}(1+O(n^{-1}))$, $g_s(d_1,\ldots,d_K)=O(n^{|w|}/N)$.
A special case is when $|w|>h$ and $|w_z|>h$ for a $z$ with $1\leq
z\leq s-1$.
In this case, no row in $\tilde H$ can match $X_s$ and $g_s(d_1,\ldots,d_K)=0$.

Thus
\[
g_s(d_1,\ldots,d_K) = \cases{
1+O\bigl(n^{-1}\bigr), &\quad $|w|< h,$
\vspace*{2pt}\cr
O(1), &\quad $|w|= h,$
\vspace*{2pt}\cr
0, & \quad$|w|>h, \max\bigl(|w_1|,\ldots,|w_{s-1}|\bigr)>h,$
\vspace*{2pt}\cr
O\bigl(n^{|w|}/N\bigr), &\quad $\mbox{otherwise.}$}
\]
Furthermore, 
the value of $g_s(d_1,\ldots,d_K)$ is the same in any regions defined
by $D_{i_1}^1 \times\cdots\times D_{i_K}^K$ in which $i_k=0,1,\ldots,s-1$
for $k=1,\ldots,K$.
Thus, write
\[
g_s(d_1,\ldots,d_K) = \sum
_{i_1,\ldots,i_K=0}^{s-1} b_s(i_1,
\ldots,i_K,M_{s-1}) I\bigl(d_1\in
D_{i_1}^1,\ldots,d_K\in D_{i_K}^K
\bigr),
\]
where
\[
b_s(i_1,\ldots,i_K,M_{s-1}) = \cases{
1+O\bigl(n^{-1}\bigr), &\quad$ |w|< h,$
\vspace*{2pt}\cr
O(1), & \quad$|w|= h,$
\vspace*{2pt}\cr
0, &\quad $|w|>h, \max\bigl(|w_1|,\ldots,|w_{s-1}|\bigr) >h,$
\vspace*{2pt}\cr
O\bigl(n^{|w|}/N\bigr), &\quad $\mbox{otherwise,}$}
\]
and $b_s(\cdot)$ is a deterministic function on $d_1,\ldots,d_K$ and $M_{s-1}$.

\subsection{Proof of Lemma~\texorpdfstring{\protect\ref{lem:ROA:r}}{3.4}} \label{sec:proof:ROA:r}

The idea to prove Lemma~\ref{lem:ROA:r} is as follows. Note that
%
\begin{equation}
E_\mathrm{ ROA} \bigl\{\bigl(N^{1/2} \bar R\bigr)^p\bigr\}
= N^{-p/2} \sum_{a_1+\cdots+a_N=p, a_1,\ldots,a_N \geq0} E_\mathrm{
ROA}
\Biggl(\prod_{i=1}^{N} r_i^{a_i}
\Biggr). \label{eqn:EROA}
\end{equation}

Let $t$ be the number of $a_i$'s being one and $s$ be the number of
nonzero $a_i$'s;
there are at most $O(N^s)$ terms in (\ref{eqn:EROA}).
Thus it suffices to show that for any $s \leq p$,
\[
E_\mathrm{ ROA} \Biggl( \prod_{i=1}^{s}
r_i^{a_i} \Biggr) - E_\mathrm{ IID} \Biggl( \prod
_{i=1}^{s} r_i^{a_i} \Biggr) = o
\bigl(N^{p/2-s}\bigr).
\]

If $t=0$, then $s \leq p/2$. From Lemma~\ref{lem:ROA:f},
\[
E_\mathrm{ ROA} \Biggl( \prod_{i=1}^{s}
r_i^{a_i} \Biggr) - E_\mathrm{ IID} \Biggl( \prod
_{i=1}^{s} r_i^{a_i} \Biggr) = O
\bigl(n^{-1}\bigr) = o\bigl(N^{p/2-s}\bigr).
\]
If $t>0$, $E_\mathrm{ IID}  (\prod_{i=1}^{s} r_i^{a_i}  )=0$.
Thus it suffices to show that for any $1 \leq t \leq s \leq p$,
$t+a_{t+1}+\cdots+a_s=p, a_{t+1},\ldots,a_s > 1$,
\[
E_\mathrm{ ROA} \Biggl(\prod_{i=1}^{t}
r_i \prod_{i=t+1}^{s}
r_i^{a_i} \Biggr) = o\bigl(N^{p/2-s}\bigr).
\]

Because $t + 2(s-t) \leq p$, $-t/2 \leq p/2-s$.
Since $r_i = \sum_{|u|>h} f_u(x_i)$, and we can rearrange the order of
$\prod_{i=1}^t r_i$ by sorting $|u_i|$, it suffices to show for any $1
\leq t \leq s $, $|u_1|\geq|u_2| \geq\cdots\geq|u_t|>h$ and
continuous functions $f, q_{t+1}, \ldots, q_{s}$,
\[
E_\mathrm{ ROA} \Biggl\{ \prod_{i=1}^{t}
f_{u_i}(x_i) \prod_{i=t+1}^{s}
q_i(x_i) \Biggr\} = o\bigl(N^{-t/2}\bigr).
\label{eqn:EROAmEIIDold}
\]

From Lemma~\ref{lem:ROA:f}, if $s>t$,
\begin{eqnarray*}
&& E_\mathrm{ ROA} \Biggl\{ \prod_{i=1}^{t}
f_{u_i}(x_i) \prod_{i=t+1}^{s}
q_i(x_i) \Biggr\}
\\
&&\qquad= E_\mathrm{ ROA} \Biggl\{ \prod_{i=1}^{t}
f_{u_i}(x_i) \prod_{i=t+1}^{s-1}
q_i(x_i) E_\mathrm{ ROA} \bigl(q_s(x_s)|x_1,
\ldots,x_{s-1}\bigr) \Biggr\}
\\
&&\qquad= E_\mathrm{ ROA} \Biggl[ \prod_{i=1}^{t}
f_{u_i}(x_i) \prod_{i=t+1}^{s-1}
q_i(x_i) \bigl\{ E_\mathrm{ IID} \bigl(q_s(x_s)
\bigr) + O\bigl(n^{-1}\bigr) \bigr\} \Biggr]
\\
&&\qquad= E_\mathrm{ ROA} \Biggl\{ \prod_{i=1}^{t}
f_{u_i}(x_i) \prod_{i=t+1}^{s-1}
q_i(x_i) \Biggr\} E_\mathrm{ IID} \bigl(q_s(x_s)
\bigr) + O\bigl(n^{-1}\bigr).
\end{eqnarray*}
Inducting on $s$, it is not hard to conclude that it is suffice to show
%
\begin{equation}
E_\mathrm{ ROA} \Biggl\{ \prod_{i=1}^{t}
f_{u_i}(x_i) \Biggr\} = o\bigl(N^{-t/2}\bigr).
\label{eqn:EROAmEIID}
\end{equation}

To show (\ref{eqn:EROAmEIID}), first express
\[
E_\mathrm{ ROA} \Biggl\{ \prod_{i=1}^{t}
f_{u_i}(x_i) \Biggr\} =  E_\mathrm{ ROA} \Biggl[ \prod
_{i=1}^{t-1} f_{u_i}(x_i)
E_\mathrm{ ROA} \bigl\{ f_{u_t}(X_t) \mid X_1,
\ldots,X_{t-1} \bigr\} \Biggr].
\]

From Lemma~\ref{lem:ROA:g},
%
\begin{eqnarray}\label{eqn:ROA:first}
&& E_\mathrm{ ROA} \Biggl\{ \prod_{i=1}^{t}
f_{u_i}(x_i) \Biggr\}
\nonumber
\\[-8pt]
\\[-8pt]
\nonumber
&&\qquad= \sum_{ i_1,\ldots,i_K }
E_\mathrm{ ROA} \Biggl\{ \prod_{i=1}^{t-1}
f_{u_i}(X_i) b_t(i_1,
\ldots,i_K,M_{t-1}) \biggl( \int_{D_t}
f_{u_t}(y) \,dy \biggr) \Biggr\},
\end{eqnarray}
where $D_t=D_{i_1}^1 \times\cdots\times D_{i_K}^K$ and $D_i^k=\delta
_n(X_i^k)$.
From (\ref{eqn:int:0}),
\[
\int_{\tilde D^1\times\cdots\times\tilde D^K} f_{u}(y) \,dy=0
\]
if there is at least one $k$ such that $\tilde D^k=[0,1)$ and $k\in u$.
Therefore,
\[
\int_{D_0^1 \times\tilde D^2 \times\cdots\times\tilde D^K} f_{u_t}(y) \,dy = - \sum
_{j=1}^{t-1} \int_{\delta_n(X_j^1) \times\tilde
D^2 \times\cdots\times\tilde D^K}
f_{u_t}(y) \,dy.
\]
Consequently, $\int_{D_t} f_{u_t}(y) \,dy$ has order $O(n^{-|w\cup
u_t|})$ where $w(d_1,\ldots,d_K)= \{k\dvtx d_k>0\}$, and
(\ref{eqn:ROA:first}) has order $O(N^{-1})$.

We can further reduce the order of (\ref{eqn:ROA:first}) if $t>1$.
For any term in the sum of (\ref{eqn:ROA:first}),
%
\begin{eqnarray}\label{eqn:ROA:second}
&& E_\mathrm{ ROA} \Biggl\{ \prod_{i=1}^{t-1}
f_{u_i}(X_i) b_t(i_1,
\ldots,i_K,M_{t-1}) \biggl( \int_{D_t}
f_{u_t}(y) \,dy \biggr) \Biggr\}
\nonumber
\\
&&\qquad= \sum_{ j_1,\ldots,j_K } E_\mathrm{ ROA} \Biggl[ \prod
_{i=1}^{t-2} f_{u_i}(X_i)
b_{t-1}(j_1,\ldots,j_K,M_{t-2})
\nonumber
\\[-8pt]
\\[-8pt]
\nonumber
&&\hspace*{87pt}\qquad{}\times \biggl\{ \int_{D_{t-1}} b_t(i_1,
\ldots,i_K,M_{t-1}) \\
&&\hspace*{150pt}{}\times\biggl( \int_{D_t}
f_{u_t}(y_t) \,dy_t \biggr) f_{u_{t-1}}(y_{t-1})
\,dy_{t-1} \biggr\} \Biggr], \nonumber
\end{eqnarray}
where $D_{t-1}=D_{j_1}^1 \times\cdots\times D_{j_K}^K$ and
$D_j^k=\delta_n(X_i^k)$.
In any region $D_{t-1}$, $b_t(\cdot)$ becomes a deterministic function on
$i_1,\ldots,i_K,M_{t-2}$ with the same order as in Lemma~\ref{lem:ROA:g}.
Let $b_t^\prime(\cdot)$ denote this function.
Then
\begin{eqnarray*}
&& E_\mathrm{ ROA} \Biggl\{ \prod_{i=1}^{t-1}
f_{u_i}(X_i) b_t(i_1,
\ldots,i_K,M_{t-1}) \biggl( \int_{D_t}
f_{u_t}(y) \,dy \biggr) \Biggr\}
\\
&&\qquad= \sum_{ j_1,\ldots,j_K } E_\mathrm{ ROA} \Biggl[ \prod
_{i=1}^{t-2} f_{u_i}(X_i)
b_{t-1}(j_1,\ldots,j_K,M_{t-2})
b_t^\prime(i_1,\ldots,i_K,M_{t-2})
\\
&&\hspace*{113pt}\qquad{}\times \biggl\{ \int_{D_{t-1}} \biggl( \int_{D_t}
f_{u_t}(y_t) \,dy_t \biggr) f_{u_{t-1}}(y_{t-1})
\,dy_{t-1} \biggr\} \Biggr].
\end{eqnarray*}

So far we have showed the first two steps to reduce the order of
magnitudes for $E_\mathrm{ ROA}\{ \prod_{i=1}^{t} f_{u_i}(x_i) \}$.
In (\ref{eqn:ROA:first}), we took $f_{u_t}(X_t)$ out of the product and
reached the $O(N^{-1})$ order.
We keep taking out\vadjust{\goodbreak} the $f_{u_i}(X_i)$ terms as in (\ref
{eqn:ROA:second}) and work on a more general formula as follows:
%
\begin{eqnarray}\label{eqn:ROA:G}
&&\Biggl(\prod_{l=1}^L |D_l|
\Biggr) ^{-1} E_\mathrm{ ROA} \Biggl[ \prod
_{i
=1}^t f_{u_i}(X_i)
\rho(M_t)
\nonumber
\\[-8pt]
\\[-8pt]
\nonumber
&&\hspace*{82pt}\qquad{}\times \Biggl\{ \int_{ \prod_{l=1}^L D_l } \Biggl( \prod
_{l=1}^L f_{v_l}(y_l)
\Biggr) \,dy_1 \cdots \,dy_L \Biggr\} \Biggr],
\end{eqnarray}
where $\rho(M_t)$ is a deterministic function on $M_t$ which has order
$O(1)$ for any $M_t$.
Suppose $G$ is an arbitrary term by (\ref{eqn:ROA:G}) with the
following parameters:
$0\leq t\leq p$,
$|u_1| \geq|u_2| \geq\cdots\geq|u_t| >h$, $L$ is a nonnegative integer,
$v_l \subseteq\{1,\ldots,K\}$, $|v_l|>h$,
$D_l=D_l^1 \times\cdots\times D_l^K$
and $D_l^k$ is either $[0,1)$ or $\delta_n(X_i^k)$ with $1 \leq i \leq
t$, or $\delta_n(y_i^k)$ with $l < i \leq L$.
Suppose that $C$ is an $t\times K$ zero--one matrix with the $(i,k)$th
element being one if and only if $k\in u_i$ and for any $1\leq l \leq
L$, $D_l^k \neq\delta_n(X_i^k)$.
Let $c_i$ be the number of ones in the $i$th row of $C$, and let $\theta
=\sum_{i=1}^t c_i/|u_i|$.
The following two lemmas give the orders of $G$ by the number of ones
in $C$.

\begin{lemma}\label{lem:ROA:G}
The quantity $G$ has order $O(N^{-\theta/2})$.
\end{lemma}

\begin{pf}
We show this by induction on $t$.
If $t=0$, then $\theta=0$, and the result clearly holds.
Next, assume the result holds for $t=0,\ldots,z-1$ with $z\geq1$.
It suffices to show the result holds for $t=z$.
Express
\begin{eqnarray*}
G &=& \Biggl(\prod_{l=1}^L
|D_l| \Biggr) ^{-1} E_\mathrm{ ROA} \Biggl[ \prod
_{i=1}^t f_{u_i}(X_i)
\rho(M_t) \Biggl\{ \int_{ \prod_{l=1}^L D_l } \prod
_{l=1}^L f_{v_l}(y_l)
\,dy_1 \cdots \,dy_L \Biggr\} \Biggr]
\\
&=& \Biggl(\prod_{l=1}^L
|D_l| \Biggr) ^{-1} \\
&&{}\times E_\mathrm{ ROA} \Biggl[ \prod
_{i=1}^{t-1} f_{u_i}(X_i)
\\
&&\hspace*{58pt}{}\times E_\mathrm{ ROA} \Biggl\{ \rho(M_t) f_{u_t}(X_t)
\Biggl( \int_{ \prod_{l=1}^L D_l } \prod_{l=1}^L
f_{v_l}(y_l) \,dy_1 \cdots \,dy_L
\Biggr)\Big \mid\\
&&\hspace*{248pt}{}\{ X_1,\ldots,X_{t-1} \} \Biggr\} \Biggr].
\end{eqnarray*}

From Lemma~\ref{lem:ROA:g} and similar to (\ref{eqn:ROA:first}) and
(\ref{eqn:ROA:second}),
\begin{eqnarray*}
&& E_\mathrm{ ROA} \Biggl\{ \rho(M_t) f_{u_t}(X_t)
\Biggl( \int_{ \prod
_{l=1}^L D_l } \prod_{l=1}^L
f_{v_l}(y_l) \,dy_1 \cdots \,dy_L
\Biggr) \Big\mid\{ X_1,\ldots,X_{t-1} \} \Biggr\}
\\
&&\qquad=  \int g(x_t) \rho(M_t) f_{u_t}(x_t)
\Biggl( \int_{ \prod_{l=1}^L D_l } \prod_{l=1}^L
f_{v_l}(y_l) \,dy_1 \cdots \,dy_L
\Biggr) \,dx_t
\\
&&\qquad=  \sum_{i_1,\ldots,i_K} b_t(i_1,
\ldots,i_K,M_{t-1})\\
&&\qquad\quad{}\times \int_{ D_{L+1}
}
\rho(M_t) f_{u_t}(x_t) \Biggl( \int
_{ \prod_{l=1}^L D_l } \prod_{l=1}^L
f_{v_l}(y_l) \,dy_1 \cdots \,dy_L
\Biggr) \,dx_t,
\end{eqnarray*}
where $g(x_t)$ is the conditional density of $X_t$,
$D_{L+1} = D_{L+1}^1 \times\cdots\times D_{L+1}^K$,
$D_{L+1}^k=\delta_n(X_{i_k}^k)$ if $i_k>0$
and $D_{L+1}^k=[0,1) \setminus\bigcup_{i=1}^{t-1} \delta_n(X_i^k)$ if $i_k=0$.

In any $D_{L+1}$, $\rho(M_t)$ is a deterministic function on $M_{t-1}$.
Let $\rho_{i_1,\ldots,i_K}(M_{t-1})$ denote this function.
Then for any $M_{t-1}$, $\rho_{i_1,\ldots,i_K}(M_{t-1})=O(1)$.
Thus
%
\begin{eqnarray}\label{eqn:ROA:G:pre}
G &=& \sum_{i_1,\ldots,i_K} \Biggl(\prod
_{l=1}^L |D_l| \Biggr)
^{-1} \nonumber\\
&&{}\times E_\mathrm{ ROA} \Biggl[ \prod_{i=1}^{t-1}
f_{u_i}(X_i) b_t(i_1,
\ldots,i_K,M_{t-1}) \rho_{i_1,\ldots,i_K}(M_{t-1})
\\
&&\hspace*{58pt}{}\times \int_{ D_{L+1} } f_{u_t}(x_t) \Biggl(
\int_{ \prod_{l=1}^L D_l } \prod_{l=1}^L
f_{v_l}(y_l) \,dy_1 \cdots \,dy_L
\Biggr) \,dx_t \Biggr].\nonumber
\end{eqnarray}

If $i_1=0$, $D_{L+1}^1=[0,1) \setminus\bigcup_{i=1}^{t-1} \delta_n(X_i^1)$.
If additionally $k\notin u_t$,
\begin{eqnarray*}
&& \int_{ D_{L+1} } f_{u_t}(x_t) \Biggl(
\int_{ \prod_{l=1}^L D_l } \prod_{l=1}^L
f_{v_l}(y_l) \,dy_1 \cdots \,dy_L
\Biggr) \,dx_t
\\
&&\qquad= \rho_{0}(M_{t-1}) \\
&&\qquad\quad{}\times\int_{ [0,1) \times D_{L+1}^2 \times\cdots
\times D_{L+1}^K }
f_{u_t}(x_t) \Biggl( \int_{ \prod_{l=1}^L D_l } \prod
_{l=1}^L f_{v_l}(y_l)
\,dy_1 \cdots \,dy_L \Biggr) \,dx_t,
\end{eqnarray*}
where $\rho_{0}(M_{t-1})=O(1)$.
If $k \in u_t$,
\begin{eqnarray*}
&& \int_{ D_{L+1} } f_{u_t}(x_t) \Biggl(
\int_{ \prod_{l=1}^L D_l } \prod_{l=1}^L
f_{v_l}(y_l) \,dy_1 \cdots \,dy_L
\Biggr) \,dx_t
\\
&&\qquad= \int_{ [0,1) \times D_{L+1}^2 \times\cdots\times D_{L+1}^K } f_{u_t}(x_t) \Biggl(
\int_{ \prod_{l=1}^L D_l } \prod_{l=1}^L
f_{v_l}(y_l) \,dy_1 \cdots \,dy_L
\Biggr) \,dx_t
\\
&&\qquad\quad{} - \sum_{j=1}^{t-1} \Biggl\{
\rho_j(M_{t-1})\\
&&\hspace*{66pt}{}\times \int_{ D_j^1 \times
D_{L+1}^2 \times\cdots\times D_{L+1}^K }
f_{u_t}(x_t) \\
&&\hspace*{162pt}{}\times\Biggl( \int_{ \prod_{l=1}^L D_l } \prod
_{l=1}^L f_{v_l}(y_l)
\,dy_1 \cdots \,dy_L \Biggr) \,dx_t \Biggr\},
\end{eqnarray*}
where $\rho_{j}(M_{t-1})=O(1)$ for $j=1,\ldots,t-1$.
Let
\[
\tilde b_t(i_1,\ldots,i_k) = \cases{ %
1, &\quad $|w|\leq h,$
\vspace*{2pt}\cr
n^{|w|}/N, &\quad $|w|>h.$ }
\]
Then $\tilde b_t$ is not related to $M_{t-1}$ and
$\tilde b_t(0,i_2,\ldots,i_K) \leq\tilde b_t(j,i_2,\ldots,i_K)$ for
any $j>0$.

From the arguments above, it suffices to show
%
\begin{eqnarray}\label{eqn:ROA:G:lr}
&&J(i_1,\ldots,i_K)\nonumber\\
&&\qquad= \tilde b_t(i_1,
\ldots,i_K) \Biggl(\prod_{l=1}^L
|D_l| \Biggr) ^{-1}
\nonumber
\\[-8pt]
\\[-8pt]
\nonumber
&&\qquad\quad{}\times E_\mathrm{ ROA} \Biggl\{ \prod
_{i=1}^{t-1} f_{u_i}(X_i)
\rho(M_{t-1})
\\
&&\hspace*{41pt}\qquad\quad{}\times \int_{ D_{L+1} } f_{u_t}(x_t) \Biggl(
\int_{ \prod_{l=1}^L D_l } \prod_{l=1}^L
f_{v_l}(y_l) \,dy_1 \cdots \,dy_L
\Biggr) \,dx_t \Biggr\}\nonumber
\end{eqnarray}
has order $O(N^{-\theta/2})$ for any $i_1,\ldots,i_K=0,1,\ldots,t-1$,
$\rho(M_{t-1})=O(1)$,\break
$D_{L+1} = D_{L+1}^1 \times\cdots\times D_{L+1}^K$,
$D_{L+1}^k=\delta_n(X_{i_k}^k)$ if $i_k>0$, $k=1,\ldots,K$,
$D_{L+1}^1=[0,1)$ if $i_1=0$
and $D_{L+1}^k=[0,1) \setminus\bigcup_{j=1}^{t-1} \delta_n(X_j^k) $ if
$i_k=0$, $k=2,\ldots,K$.

From similar arguments, it suffices to show (\ref{eqn:ROA:G:lr}) has
order $O(N^{-\theta/2})$ for any $i_1,\ldots,i_K=0,1,\ldots,t-1$, $\rho
(M_{t-1})=O(1)$,
$D_{L+1} = D_{L+1}^1 \times\cdots\times D_{L+1}^K$,
$D_{L+1}^k=\delta_n(X_{i_k}^k)$ if $i_k>0$ and
$D_{L+1}^k=[0,1)$ if $i_k=0$, $k=1,\ldots,K$.

If $i_k \neq0$ and $k \notin u_t$,
then any term that can be written as $J(i_1,
\ldots,i_K)$
has smaller or the same order than a term that can be written as
$J(i_1,\ldots,i_{k-1},0,\break i_{k+1}, \ldots,i_K)$.
If $i_k=0$ and the $(t,k)$th element of $C$ is one, from (\ref
{eqn:int:0}), $J=0$.
Thus it suffices to consider $J(i_1,\ldots,i_K)$ with
$i_k=0,\ldots,t-1$ for $k \in u_t$ and the $(t,k)$th element of $C$
being zero,
$i_k=1,\ldots,t-1$ for $k \in u_t$ and the $(t,k)$th element of $C$
being one
and $i_k=0$ for $k \notin u_t$.
Clearly, $w \subseteq u_t$ and $c_t \leq|w| \leq|u_t|$.

Let
\begin{eqnarray*}
G^\prime_{i_1,\ldots,i_K} &=& \Biggl(\prod_{l=1}^{L+1}
|D_l| \Biggr) ^{-1} E_\mathrm{ ROA} \Biggl[ \prod
_{i=1}^{t-1} f_{u_i}(X_i)
\rho(M_{t-1})
\\
& &\hspace*{89pt}{}\times\Biggl\{ \int_{ (\prod_{l=1}^L D_l^\prime) \times D_{L+1}
} \prod_{l=1}^{L+1}
f_{v_l}(y_l) \,dy_1 \cdots \,dy_{L+1}
\Biggr\} \Biggr],
\end{eqnarray*}
where $v_{L+1}=u_t$ and
\[
D_l^{k\prime} = %
\cases{ \delta_n
\bigl(y_{L+1}^k\bigr), &\quad $\mbox{if } D_l^k=
\delta_n\bigl(X_t^k\bigr)$,
\cr
D_l^k, &\quad  $\mbox{otherwise}$.} %
\]
Then $J$ in (\ref{eqn:ROA:G:lr}) can be expressed as
\[
J(i_1,\ldots,i_K)=\tilde b_t(i_1,
\ldots,i_K) n^{-|w|} G^\prime _{i_1,\ldots,i_K}.
\]
For any $(i_1,\ldots,i_K)$, $G^\prime_{i_1,\ldots,i_K}$ is a term by
(\ref{eqn:ROA:G}).
Furthermore, the matrix associated with $G^\prime_{i_1,\ldots,i_K}$,
denoted as $C^\prime_{i_1,\ldots,i_K}$, is a $(t-1)\times K$ matrix
with equal or fewer elements of ones than the first $t-1$ rows of $C$.
If $i_k=z>0$, the $(z,k)$th element of $C^\prime_{i_1,\ldots,i_K}$ is zero.
Other elements of $C^\prime(D_{L+1})$ are the same with that of $C$.
Let $c_i^\prime$ be the number of ones in the $i$th row of $C^\prime
_{i_1,\ldots,i_K}$,
and let $\theta^\prime= \sum_{i=1}^{t-1} c_i^\prime/ |u_i|$, and
we have
\[
\theta^\prime\geq\theta- c_t / |u_t| - |w| /
|u_t|.
\]

By induction,
%
\begin{equation}
G^\prime_{i_1,\ldots,i_K} = O\bigl(N^{- ( \theta-c_t/|u_t|-|w|/|u_t| ) /2}\bigr) =O
\bigl(N^{-\theta/2+1}\bigr) \label{eqn:ROA:G:r}
\end{equation}
and
%
\begin{equation}
G^\prime_{i_1,\ldots,i_K} = O\bigl(N^{- ( \theta-c_t/|u_t|-|w|/|u_t| ) /2}\bigr) =O
\bigl(N^{-\theta/2}n^{|w|}\bigr). \label{eqn:ROA:G:l}
\end{equation}
Consequently, $J$ in (\ref{eqn:ROA:G:lr}) has order $O(N^{-\theta/2})$.
This completes the proof.
\end{pf}

The result of Lemma~\ref{lem:ROA:G} is improved by Lemma~\ref{lem:ROA:Gi}.

\begin{lemma}\label{lem:ROA:Gi}
If $c_t>0$, $G$ has order $o(N^{-\theta/2})$.
\end{lemma}

\begin{pf}
It suffices to show $J$ in (\ref{eqn:ROA:G:lr}) has order $o(N^{-\theta/2})$.
Since $c_t>0$, (\ref{eqn:ROA:G:l})~becomes
\[
G^\prime_{i_1,\ldots,i_K} 
= O\bigl(N^{-\theta/2 + |w|/(h+1)}
\bigr) =o\bigl(N^{-\theta/2}n^{|w|}\bigr). \label{eqn:ROA:G:l:new}
\]
Therefore, for $|w|\leq h$, $J=o(N^{-\theta/2})$.
When $|w|>h$ and $c_t<|u_t|$, (\ref{eqn:ROA:G:r}) becomes
%
\begin{equation}
G^\prime_{i_1,\ldots,i_K} = O\bigl(N^{- ( \theta-c_t/|u_t|-|w|/|u_t| ) /2}\bigr) =o
\bigl(N^{-\theta/2+1}\bigr), \label{eqn:ROA:G:l:new2}
\end{equation}
and $J=o(N^{-\theta/2})$.
When $|w|>h$ and there is a $j$ such that $|w_j|>h$, $b_t(i_1,\ldots,i_K)$ in (\ref{eqn:ROA:G:pre}) is zero and $J=0$.

It remains to show $G^\prime_{i_1,\ldots,i_K}=o(N^{-\theta/2+1})$ for
$c_t=|w|=|u_t|$ and\break  $\max(|w_i|)\leq h$.
Let $\{(j_1,k_1),\ldots,(j_z,k_z)\}$ denote the elements of $C^\prime
_{d_1,\ldots,d_K}$ that are different from those of the first $t-1$
rows of $C$.
When $|u_{j_x}|>|u_t|$ for an $x$ with $1\leq x\leq z$, (\ref
{eqn:ROA:G:r}) becomes
\[
G^\prime_{d_1,\ldots,d_K}=O\bigl(N^{- \{ \theta
-1/|u_{j_x}|-(c_t+|w|-1)/|u_t| \} /2}\bigr) =o
\bigl(N^{-\theta/2+1}\bigr).
\]
When $z<|w|$, (\ref{eqn:ROA:G:r}) becomes
\[
G^\prime_{i_1,\ldots,i_K}=O\bigl(N^{- ( \theta-z/|u_t|-|w|/|u_t| ) /2}\bigr) =o
\bigl(N^{-\theta/2+1}\bigr).
\]
Finally, when $z=|w|$ and $|u_{j_1}|=\cdots=|u_{j_z}|=|u_t|$,
since $\max_j\{|w_j|\}\leq h$, $\{j_1,\ldots,j_z\}$ are not all equal
to each other.
Consequently, there is at least one $x$ such that $0< c_{j_x}^\prime<
|u_{j_x}|=|u_t|$.
From (\ref{eqn:ROA:G:l:new2}), $G^\prime_{d_1,\ldots,d_K}=o(N^{-\theta/2+1})$.
Combining all cases, $J=o(N^{-\theta/2})$.
This completes the proof.
\end{pf}

We now give the proof of Lemma~\ref{lem:ROA:r}.

\begin{pf}
We have argued in (\ref{eqn:EROAmEIID}) that it suffices to show for
any $1 \leq t \leq p $, $|u_1|\geq|u_2| \geq\cdots\geq|u_t|>h$ and
continuous functions $f$,
\[
E_\mathrm{ ROA} \Biggl\{ \prod_{i=1}^{t}
f_{u_i}(x_i) \Biggr\} = o\bigl(N^{-t/2}\bigr).
\]

$ E_\mathrm{ ROA} \{\prod_{i=1}^{t} f_{u_i}(x_i)  \} $ is a term
by (\ref{eqn:ROA:G}) with $\theta=t$ and $c_t=|u_t|>0$.
From Lemma~\ref{lem:ROA:Gi}, $E_\mathrm{ ROA}  \{ \prod_{i=1}^{t}
f_{u_i}(x_i)  \} =o(N^{-t/2}) $.
This completes the proof.
\end{pf}

\subsection{Proof of Lemma~\texorpdfstring{\protect\ref{lem:UD:g}}{4.1}} \label{sec:proof:UD:g}

Similar to the argument in the proof of Lem\-ma~\ref{lem:ROA:g}, we have that
\begin{eqnarray*}
&&g_s(d_1,
\ldots,d_K)\\
&&\qquad
= \cases{ %
1+O\bigl(n^{-1}\bigr), &\quad $|w|< h,$
\vspace*{2pt}\cr
O(1), &\quad $|w|= h,$
\vspace*{2pt}\cr
0, &\quad $|w|>h, \max\bigl(|w_1|,\ldots,|w_{s-1}|\bigr)>h,$
\vspace*{2pt}\cr
O\bigl(n^{|w|}/N\bigr), & \quad$\mbox{otherwise.}$}
\end{eqnarray*}

However, a special case is when there is a $k$ such that $i_k>s-1$.
From~(\ref{eqn:const:UD}), two rows cannot be in the same subdivision
with length $1/N$.
Thus $g_s=0$ in this case.

Next, the density is uniform in each of the $D_{i_1}^1 \times\cdots
\times D_{i_K}^K$ regions,
where $i_1,\ldots,i_K = 0,\ldots,2s-2$.
Thus we can write
\[
g_s(d_1,\ldots,d_K) = \sum
_{i_1,\ldots,i_K=0}^{2s-2} b_s(i_1,
\ldots,i_K,M_{s-1}) I\bigl(d_1\in
D_{i_1}^1,\ldots,d_K\in D_{i_K}^K
\bigr),
\]
where
\[
b_s(i_1,\ldots,i_K,M_{s-1}) = \cases{
 0, &\quad $\mbox{there is a $k$ such that }
i_k>s-1,$
\vspace*{2pt}\cr
1+O\bigl(n^{-1}\bigr), &\quad $i_1,\ldots,i_K\leq
s-1, |w|< h,$
\vspace*{2pt}\cr
O(1), &\quad $i_1,\ldots,i_K\leq s-1, |w|= h,$
\vspace*{2pt}\cr
0, &\quad  $|w|>h, \max\bigl(|w_1|,\ldots,|w_{s-1}|\bigr) >h,$
\vspace*{2pt}\cr
O\bigl(n^{|w|}/N\bigr), &\quad  $\mbox{otherwise,}$}
\]
and $b_s(\cdot)$ is a deterministic function on $d_1,\ldots,d_K$ and $M_{s-1}$.

\subsection{A sketch to prove Lemma~\texorpdfstring{\protect\ref{lem:UD:r}}{4.3}} \label{sec:proof:UD:r}
Suppose $G$ is an arbitrary term given by
%
\begin{eqnarray}\label{eqn:UD:G}
&&\Biggl(\prod_{l=1}^L |D_l|
\Biggr) ^{-1}
\nonumber
\\[-8pt]
\\[-8pt]
\nonumber
&&\qquad{}\times E_\mathrm{ UD} \Biggl[ \prod
_{i
=1}^t f_{u_i}(X_i)
\rho(M_t) \Biggl\{ \int_{ \prod_{l=1}^L D_l } \Biggl( \prod
_{l=1}^L f_{v_l}(y_l)
\Biggr) \,dy_1 \cdots \,dy_L \Biggr\} \Biggr],
\end{eqnarray}
with the following parameters:
$\rho(M_t)$ is a deterministic function on $M_t$ which has order $O(1)$
for any $M_t$,
$0\leq t\leq p$,
$|u_1| \geq|u_2| \geq\cdots\geq|u_t| >h$, $L$ is a nonnegative integer,
$v_l \subseteq\{1,\ldots,K\}$, $|v_l|>h$,
$D_l=D_l^1 \times\cdots\times D_l^K$
and $D_l^k$ is either $[0,1)$ or $\delta_n(X_i^k)$ with $1 \leq i \leq
t$, or $\delta_n(y_i^k)$ with $l < i \leq L$, or $\delta_N(X_i^k)$ with
$1 \leq i \leq t$, or $\delta_N(y_i^k)$ with $l < i \leq L$.
Suppose that $C$ is an $t\times K$ zero--one matrix with the $(i,k)$th
element being one if and only if $k\in u_i$ and for any $1\leq l \leq
L$, $D_l^k$ is neither $\delta_n(X_i^k)$ nor $\delta_N(X_i^k)$.
Let $c_i$ be the number of ones in the $i$th row of $C$, and let $\theta
=\sum_{i=1}^t c_i/|u_i|$.
The following two lemmas give the order of $G$.

\begin{lemma}\label{lem:UD:G}
The quantity $G$ has order $O(N^{-\theta/2})$.
\end{lemma}

\begin{pf}
We show this by induction on $t$.
If $t=0$, then $\theta=0$, and the result clearly holds.
Next, assume the result holds for $t=0,\ldots,z-1$ with $z\geq1$.
It suffices to show the result holds for $t=z$.
Express
\begin{eqnarray*}
G &=& \Biggl(\prod_{l=1}^L
|D_l| \Biggr) ^{-1} E_\mathrm{ UD} \Biggl[ \prod
_{i=1}^t f_{u_i}(X_i)
\rho(M_t) \Biggl\{ \int_{ \prod_{l=1}^L D_l } \prod
_{l=1}^L f_{v_l}(y_l)
\,dy_1 \cdots \,dy_L \Biggr\} \Biggr]
\\
&=& \Biggl(\prod_{l=1}^L
|D_l| \Biggr) ^{-1}\\
&&{}\times E_\mathrm{ UD} \Biggl[ \prod
_{i=1}^{t-1} f_{u_i}(X_i)
E_\mathrm{ UD} \Biggl\{ \rho(M_t) f_{u_t}(X_t)\\
&&\hspace*{93pt}\qquad{}\times
\Biggl( \int_{ \prod_{l=1}^L D_l } \prod_{l=1}^L
f_{v_l}(y_l) \,dy_1 \cdots \,dy_L
\Biggr) \Big\mid\\
&&\hspace*{209pt}\{ X_1,\ldots,X_{t-1} \} \Biggr\} \Biggr]\hspace*{-1pt}.
\end{eqnarray*}

From Lemma~\ref{lem:UD:g},
\begin{eqnarray*}
&& E_\mathrm{ UD} \Biggl\{ \rho(M_t) f_{u_t}(X_t)
\Biggl( \int_{ \prod
_{l=1}^L D_l } \prod_{l=1}^L
f_{v_l}(y_l) \,dy_1 \cdots \,dy_L
\Biggr) \Big\mid\{ X_1,\ldots,X_{t-1} \} \Biggr\}
\\
&&\qquad=  \int g(x_t) \rho(M_t) f_{u_t}(x_t)
\Biggl( \int_{ \prod_{l=1}^L D_l } \prod_{l=1}^L
f_{v_l}(y_l) \,dy_1 \cdots \,dy_L
\Biggr) \,dx_t
\\
&&\qquad=  \sum_{i_1,\ldots,i_K} b_t(i_1,
\ldots,i_K,M_{t-1}) \int_{ D_{L+1}
}
\rho(M_t) f_{u_t}(x_t) \\
&&\hspace*{60pt}{}\times\Biggl( \int
_{ \prod_{l=1}^L D_l } \prod_{l=1}^L
f_{v_l}(y_l) \,dy_1 \cdots \,dy_L
\Biggr) \,dx_t,
\end{eqnarray*}
where $g(x_t)$ is the conditional density of $X_t$,
$D_{L+1} = D_{L+1}^1 \times\cdots\times D_{L+1}^K$,
$D_{L+1}^k=\delta_n(X_{i_k}^k)\setminus\bigcup_{j=1}^{t-1} \delta
_N(X_j^k)$ if $0<i_k\leq t-1$,
$D_{L+1}^k=\delta_N(X_{i_k-(t-1)}^k)$ if $i_k>t-1$
and $D_{L+1}^k=[0,1) \setminus\bigcup_{j=1}^{t-1} \delta_n(X_j^k)$ if $i_k=0$.

In any $D_{L+1}$, $\rho(M_t)$ is a deterministic function on $M_{t-1}$.
Let $\rho_{i_1,\ldots,i_K}(M_{t-1})$ denote this function.
Then for any $M_{t-1}$, $\rho_{i_1,\ldots,i_K}(M_{t-1})=O(1)$.
Thus
\begin{eqnarray*}
G &=& \sum_{i_1,\ldots,i_K} \Biggl(\prod
_{l=1}^L |D_l| \Biggr)
^{-1}\\
&&\hspace*{26pt}{}\times E_\mathrm{ UD} \Biggl\{ \prod_{i=1}^{t-1}
f_{u_i}(X_i) b_t(i_1,
\ldots,i_K,M_{t-1}) \rho_{i_1,\ldots,i_K}(M_{t-1})
\\
&&\hspace*{65pt}{}\times \int_{ D_{L+1} } f_{u_t}(x_t) \Biggl(
\int_{ \prod_{l=1}^L D_l } \prod_{l=1}^L
f_{v_l}(y_l) \,dy_1 \cdots \,dy_L
\Biggr) \,dx_t \Biggr\}.
\end{eqnarray*}

If $0<i_1\leq t-1$, $D_{L+1}^1=\delta_n(X_{i_1}^1) \setminus\bigcup_{j=1}^{t-1} \delta_N(X_j^1)$.
If additionally $k\notin u_t$,
\begin{eqnarray*}
&& \int_{ D_{L+1} } f_{u_t}(x_t) \Biggl(
\int_{ \prod_{l=1}^L D_l } \prod_{l=1}^L
f_{v_l}(y_l) \,dy_1 \cdots \,dy_L
\Biggr) \,dx_t
\\
&&\qquad= \rho_{0}^\prime(M_{t-1})\\
&&\qquad\quad{}\times \int
_{\delta_n(X_{i_1}^1) \times
D_{L+1}^2 \times\cdots\times D_{L+1}^K } f_{u_t}(x_t) \\
&&\hspace*{145pt}{}\times\Biggl( \int
_{ \prod_{l=1}^L D_l } \prod_{l=1}^L
f_{v_l}(y_l) \,dy_1 \cdots \,dy_L
\Biggr) \,dx_t,
\end{eqnarray*}
where $\rho_{0}^\prime(M_{t-1})=O(1)$.
If $k \in u_t$,
\begin{eqnarray*}
&& \int_{ D_{L+1} } f_{u_t}(x_t) \Biggl(
\int_{ \prod_{l=1}^L D_l } \prod_{l=1}^L
f_{v_l}(y_l) \,dy_1 \cdots \,dy_L
\Biggr) \,dx_t
\\
&&\qquad= \int_{ \delta_n(X_{i_1}^1) \times D_{L+1}^2 \times\cdots\times
D_{L+1}^K } f_{u_t}(x_t) \Biggl(
\int_{ \prod_{l=1}^L D_l } \prod_{l=1}^L
f_{v_l}(y_l) \,dy_1 \cdots \,dy_L
\Biggr) \,dx_t
\\
&&\qquad\quad{} - \sum_{j=1}^{t-1} \Biggl\{ \int
_{ \delta_N({X_j^1}) \times D_{L+1}^2
\times\cdots\times D_{L+1}^K } f_{u_t}(x_t) \\
&&\hspace*{167pt}{}\times\Biggl( \int
_{ \prod_{l=1}^L D_l } \prod_{l=1}^L
f_{v_l}(y_l) \,dy_1 \cdots \,dy_L
\Biggr) \,dx_t \Biggr\}.\
\end{eqnarray*}

If $i_1=0$, $D_{L+1}^1=[0,1) \setminus\bigcup_{i=1}^{t-1} \delta_n(X_i^1)$.
If additionally $k\notin u_t$,
\begin{eqnarray*}
&& \int_{ D_{L+1} } f_{u_t}(x_t) \Biggl(
\int_{ \prod_{l=1}^L D_l } \prod_{l=1}^L
f_{v_l}(y_l) \,dy_1 \cdots \,dy_L
\Biggr) \,dx_t
\\
&&\qquad= \rho_{0}^\prime(M_{t-1})\\
&&\qquad\quad{}\times \int
_{ [0,1) \times D_{L+1}^2 \times
\cdots\times D_{L+1}^K } f_{u_t}(x_t) \Biggl( \int
_{ \prod_{l=1}^L D_l } \prod_{l=1}^L
f_{v_l}(y_l) \,dy_1 \cdots \,dy_L
\Biggr) \,dx_t,
\end{eqnarray*}
where $\rho_{0}^\prime(M_{t-1})=O(1)$.
If $k \in u_t$,
\begin{eqnarray*}
\hspace*{-4pt}&& \int_{ D_{L+1} } f_{u_t}(x_t) \Biggl(
\int_{ \prod_{l=1}^L D_l } \prod_{l=1}^L
f_{v_l}(y_l) \,dy_1 \cdots \,dy_L
\Biggr) \,dx_t
\\
\hspace*{-4pt}&&\qquad= \int_{ [0,1) \times D_{L+1}^2 \times\cdots\times D_{L+1}^K } f_{u_t}(x_t) \Biggl(
\int_{ \prod_{l=1}^L D_l } \prod_{l=1}^L
f_{v_l}(y_l) \,dy_1 \cdots \,dy_L
\Biggr) \,dx_t
\\
\hspace*{-4pt}&&\qquad\quad{} - \sum_{j=1}^{t-1} \Biggl\{
\rho_j^\prime(M_{t-1})\\
\hspace*{-4pt}&&\hspace*{65pt}{}\times \int_{ \delta
_n(X_j^1) \times D_{L+1}^2 \times\cdots\times D_{L+1}^K }
f_{u_t}(x_t)\\
\hspace*{-4pt}&&\hspace*{172pt}{}\times \Biggl( \int_{ \prod_{l=1}^L D_l } \prod
_{l=1}^L f_{v_l}(y_l)
\,dy_1 \cdots \,dy_L \Biggr) \,dx_t \Biggr\},
\end{eqnarray*}
where $\rho_{j}^\prime(M_{t-1})=O(1)$ for $j=1,\ldots,t-1$.
Let
\[
\tilde b_t(i_1,\ldots,i_k) = \cases{ %
1, &\quad $|w|\leq h,$
\vspace*{2pt}\cr
n^{|w|}/N, &\quad $|w|>h.$}
\]
Then $\tilde b_t$ is not related to $M_{t-1}$ and
$\tilde b_t(0,i_2,\ldots,i_K) \leq\tilde b_t(j,i_2,\ldots,i_K)$ for
any $j>0$.

From arguments above, it suffices to show
%
\begin{eqnarray}\label{eqn:UD:G:lr}
&&J(i_1,\ldots,i_K)\nonumber\\
&&\qquad= \tilde b_t(i_1,
\ldots,i_K) \Biggl(\prod_{l=1}^L
|D_l| \Biggr) ^{-1}
\nonumber
\\[-8pt]
\\[-8pt]
\nonumber
&&\qquad\quad{}\times E_\mathrm{ UD} \Biggl\{ \prod
_{i=1}^{t-1} f_{u_i}(X_i)
\rho(M_{t-1}) \\
&&\hspace*{73pt}{}\times\int_{ D_{L+1} } f_{u_t}(x_t) \Biggl(
\int_{ \prod_{l=1}^L D_l } \prod_{l=1}^L
f_{v_l}(y_l) \,dy_1 \cdots \,dy_L
\Biggr) \,dx_t \Biggr\}\nonumber
\end{eqnarray}
has order $O(N^{-\theta/2})$ for any $i_1,\ldots,i_K=0,1,\ldots,t-1$,
$\rho(M_{t-1})=O(1)$,\break
$D_{L+1} = D_{L+1}^1 \times\cdots\times D_{L+1}^K$,
$D_{L+1}^1=\delta_n(X_{i_1}^1) $ if $0<i_1\leq t-1$,
$D_{L+1}^k=\delta_n(X_{i_k}^k) \setminus\bigcup_{j=1}^{t-1} \delta
_N(X_j^k)$ if $0<i_k\leq t-1$, $k=2,\ldots,K$,
$D_{L+1}^k=\delta_N(X_{i_k-(t-1)}^k)$ if $i_k> t-1$, $k=1,\ldots,K$,
$D_{L+1}^1=[0,1)$ if $i_1=0$
and $D_{L+1}^k=[0,1) \setminus\bigcup_{j=1}^{t-1} \delta_n(X_j^k) $ if
$i_k=0$, $k=2,\ldots,K$.

From similar arguments, it suffices to show (\ref{eqn:UD:G:lr}) has
order $O(N^{-\theta/2})$ for any $i_1,\ldots,i_K=0,1,\ldots,t-1$, $\rho
(M_{t-1})=O(1)$,
$D_{L+1} = D_{L+1}^1 \times\cdots\times D_{L+1}^K$,
$D_{L+1}^k=\delta_n(X_{i_k}^k)$ if $0<i_k\leq t-1$,
$D_{L+1}^k=\delta_N(X_{i_k-(t-1)}^k)$ if $i_k>t-1$ and
$D_{L+1}^k=[0,1)$ if $i_k=0$, $k=1,\ldots,K$.

If $i_k \neq0$ and $k \notin u_t$,
then any term that can be written as $J(i_1,
\ldots,i_K)$
has smaller or the same order than a term that can be written as
$J(i_1,\ldots,i_{k-1},0,\break i_{k+1}, \ldots,i_K)$.
If $i_k=0$ and the $(t,k)$th element of $C$ is one, from (\ref
{eqn:int:0}), $J=0$.
Thus it suffices to consider $J(i_1,\ldots,i_K)$ with
$i_k=0,\ldots,2t-2$ for $k \in u_t$ and the $(t,k)$th element of $C$
being zero,
$i_k=1,\ldots,2t-2$ for $k \in u_t$ and the $(t,k)$th element of $C$
being one
and $i_k=0$ for $k \notin u_t$.
Clearly, $w \subseteq u_t$ and $c_t \leq|w| \leq|u_t|$.

Let
\begin{eqnarray*}
G^\prime_{i_1,\ldots,i_K} &=& \Biggl(\prod_{l=1}^{L+1}
|D_l| \Biggr) ^{-1} E_\mathrm{ UD} \Biggl[ \prod
_{i=1}^{t-1} f_{u_i}(X_i)
\rho(M_{t-1})
\\
&&\hspace*{86pt}{}\times \Biggl\{ \int_{ (\prod_{l=1}^L D_l^\prime) \times D_{L+1}
} \prod_{l=1}^{L+1}
f_{v_l}(y_l) \,dy_1 \cdots \,dy_{L+1}
\Biggr\} \Biggr],
\end{eqnarray*}
where $v_{L+1}=u_t$ and
\[
D_l^{k\prime} = %
\cases{ \delta_n
\bigl(y_{L+1}^k\bigr), &\quad $\mbox{if } D_l^k=
\delta_n\bigl(X_t^k\bigr)$,
\cr
\delta_N\bigl(y_{L+1}^k\bigr), &\quad $\mbox{if }
D_l^k= \delta_N\bigl(X_t^k
\bigr)$,
\cr
D_l^k, &\quad $\mbox{otherwise}$.} %
\]
Then $J$ in (\ref{eqn:UD:G:lr}) can be expressed as
\[
J(i_1,\ldots,i_K)=\tilde b_t(i_1,
\ldots,i_K) n^{-|w|} G^\prime _{i_1,\ldots,i_K}.
\]
For any $(i_1,\ldots,i_K)$, $G^\prime_{i_1,\ldots,i_K}$ is a term by
(\ref{eqn:UD:G}).
Furthermore, the matrix associated with $G^\prime_{i_1,\ldots,i_K}$,
denoted as $C^\prime_{i_1,\ldots,i_K}$, is a $(t-1)\times K$ matrix
with equal or fewer elements of ones than the first $t-1$ rows of $C$.
If $0<i_k=z\leq t-1$, the $(z,k)$th element of $C^\prime_{i_1,\ldots,i_K}$ is zero.
If $i_k=z> t-1$, the $(z-(t-1),k)$th element of $C^\prime_{i_1,\ldots,i_K}$ is zero.
Other elements of $C^\prime(D_{L+1})$ are the same with that of $C$.
Let $c_i^\prime$ be the number of ones in the $i$th row of $C^\prime
_{i_1,\ldots,i_K}$,
and let $\theta^\prime= \sum_{i=1}^{t-1} c_i^\prime/ |u_i|$, so
we have
\[
\theta^\prime\geq\theta- c_t / |u_t| - |w| /
|u_t|.
\]

By induction,
%
\begin{equation}
G^\prime_{i_1,\ldots,i_K} = O\bigl(N^{- ( \theta-c_t/|u_t|-|w|/|u_t| ) /2}\bigr) =O
\bigl(N^{-\theta/2+1}\bigr) \label{eqn:UD:G:r}
\end{equation}
and
%
\begin{equation}
G^\prime_{i_1,\ldots,i_K} = O\bigl(N^{- ( \theta-c_t/|u_t|-|w|/|u_t| ) /2}\bigr) =O
\bigl(N^{-\theta/2}n^{|w|}\bigr). \label{eqn:UD:G:l}
\end{equation}
Consequently, $J$ in (\ref{eqn:UD:G:lr}) has order $O(N^{-\theta/2})$.
This completes the proof.
\end{pf}

\begin{lemma}\label{lem:UD:Gi}
If $c_t>0$, $G$ has order $o(N^{-\theta/2})$.
\end{lemma}

The proofs for Lemma~\ref{lem:UD:Gi} and~\ref{lem:UD:r} are similar to
the proofs for Lemma~\ref{lem:ROA:Gi} and~\ref{lem:ROA:r},
respectively, and are omitted.
\end{appendix}

\section*{Acknowledgments}
The authors thank the Editor, the Associate Editors and two referees
for their valuable comments and suggestions
that improved this article.
%


%



\printaddresses
\end{document}